\documentclass[12pt]{article}

\usepackage{amssymb,amsmath,amsfonts,amsthm}
\usepackage[dvips]{graphicx}
\usepackage{subfigure}

\newcommand{\R}{{\mathbb R}}

\newtheorem{thh}{Theorem}
\newtheorem{lemm}{Lemma}
\newtheorem{prop}{Proposition}

\newtheorem{remarque}{Remark}

\newcommand{\ud} {\, {\mathrm{d}}}
\newcommand{\eps} {\varepsilon}
\newcommand{\e} {\mbox{e}}
\newcommand{\sign} {\mbox{sign}}

\begin{document}

\title{Convergence to the equilibrium state for \\ Bose-Einstein 1-D Kac grazing limit model}
\author{
 Radjesvarane ALEXANDRE
\thanks{ Department of Mathematics, Shanghai Jiao Tong University, Shanghai 200240, China.}
 \thanks{IRENAV, Arts et Metiers Pairs Tech, Ecole Navale, Lanveoc Poulmic, Brest 29290, France.  E-mail: {\tt radjesvarane.alexandre@ecole-navale.fr}}
\and
Jie LIAO\thanks{School of Science, East China University of Science and Technology.  E-mail:{\tt liaojie@ecust.edu.cn}}
\and
Chunjin LIN
\thanks{Department of Mathematics, Shanghai Jiao Tong University, Shanghai 200240, China. }
\thanks{ Department of Mathematics, College of Sciences, Hohai University,
Nanjing 210098, Jiangsu, China. E-mail: {\tt cjlin@hhu.edu.cn}
  }
}
\date{\today}
\maketitle

\begin{abstract}
The convergence to the equilibrium of the solution of the quantum Kac model for Bose-Einstein identical particles is studied in this paper. Using the relative entropy method and a detailed analysis of
the entropy production, the exponential decay rate is obtained under suitable assumptions.
The theoretical results are further illustrated by  numerical
simulations.
\end{abstract}

\section{Introduction}

In this paper we study the equation governing the time evolution of a gas
composed of Bose-Einstein identical particles. Let $f(t,v)$ be the velocity
 distribution function at time $t>0$ with the velocity $v\in \R.$
According to quantum physics, the presence of a particle in the velocity
range $\ud v$ increases the probability that a particle will enter that range:
the presence of $f(v)\ud v$ particles per unit volume increases this probability
in the radio $1+f(v).$ Following Chapman and Cowling \cite{chapman},
this fundamental assumption yields the so-called Boltzmann-Bose-Einstein equation,
that is the quantum Boltzmann equation for Bose-Einstein particles.
This equation has been extensively studied in physical literatures and numerical simulations. However, there are not many rigorous mathematical results. We mention here for spatial homogeneous isotropic case, a theory of weak solutions developed by Lu in \cite{lu,luzhang}, and another class of locally defined in time classical solutions by Escobedo et al. in \cite{escobedo,escobedo2,escobedo3}. See \cite{ArkNou,spo} for more reviews of currently available mathematical results. On the other hand, Allemand and Toscani in \cite{allemand} derived the following nonlinear Fokker-Planck equation (Kac model)
\begin{equation}
\label{kac}
    \partial_t f = A_f(t) \partial_{vv}f +B_f(t) \partial_v (vf(1+f))
\end{equation}
with
\[
    A_f(t)= \int v^2f(1+f) \ud v, \quad B_f(t) = \int f \ud v.
\]
This model  is obtained as the grazing collision limit of one-dimensional Boltzmann equation for Bose-Einstein particles in the spirit of Kac caricature of a Maxwell gas with a singular kernel \cite{allemand}.
However, the existence of good solutions for this integro partial differential equation is still unknown, and is currently under investigation.
We will furthermore specialize to this Kac model, and study the convergence of the solution for the Kac model to the Bose distribution by using entropy method.

The rigorous study of the convergence to equilibrium is by now classical in kinetic theory.  For example, using the classical logarithmic-Sobolev
inequality of Gross \cite{Gross},
and the Csiszar-Kullback-Pinsker inequality \cite{Csiszar,Kullback},
the convergence to the equilibrium with exponential decay rate can be
derived by the relative entropy method for linear Fokker-Plank equation. For the nonlinear Fokker-Plank-Landau equation, the trend to equilibrium has been obtained by Desvillettes and
Villani in \cite{DV2}.
Toscani and Villani in \cite{TV} studied the convergence to the equilibrium for
the Boltzmann equation.  Except for the spacial homogeneous kinetic models,
Desvillettes and Villani \cite{DV3} studied the trend to equilibrium for
the spacial inhomogeneous linear Fokker-Planck equation.
For more about the trend to equilibrium for classical kinetic equations,
we refer to \cite{DV,DV4,Des} and references therein.

In \cite{Carrillo}, Carrillo, Rosado and Salvarani studied
a 1-D quantum Fokker-planck equation
\[
    \partial_t f = \partial_{vv} f + \partial_v (vf(1+f).
\]
Note that the factor $1+f$ comes from the quantum effects.
It is easy to see that, the mass, $\int f(t,v)\ud v $,  is
conserved along the time evolution.
By using the relative entropy method, they proved that the solutions converge to the Bose equilibrium with exponential decay rate.
The above model, a simplified model of \eqref{kac} with $A_f(t)$ and $B_f(t)$ replaced by constant 1, does not conserve the kinetic energy.

For the Kac grazing limit model \eqref{kac}, in comparing with
the model studied in \cite{Carrillo}, there is an additional
conservation law: conservation of kinetic energy, i.e $\int |v|^2 f \ud v$. However, the nonlinearity in the Kac model is stronger.
For later use, let $m$ and $e$ be the mass and the energy defined by the initial data $f_0(v),$
\[
    m=\int f_0(v) \ud v, \quad e=\int v^2 f_0(v) \ud v,
\]
by supposing that $f_0>0, $ $f_0\in L^1(\R),$ and $v^2 f_0 \in L^1(\R).$
The entropy $H(f)$  is defined as
\[
    H(f)(t) =\int \gamma(f)(t,v)\ud v, \quad  \mbox{with } \gamma(x) = x\log x  - (1+x) \log(1+x).
\]
We further remark that the entropy used in \cite{Carrillo} is the sum of the entropy defined
above and the kinetic energy which is conserved for the Kac grazing limit model \eqref{kac}.

We shall work with smooth enough solutions: we show in the next section that one can get a priori weighted $L^2$ bounds on solutions, and similar estimates hold also true for higher derivatives.

The main result in this paper is stated as the following Theorem.
\begin{thh}
\label{mainthm}
Let $f(t,v)$ be the solution of the Kac grazing limit model
\eqref{kac} with initial data $f_0$  which is positive and satisfies
\[
 \int f_0(v)\ud v =m, \quad \int v^2 f_0(v) \ud v =e
\]
for some positive constants $m,\ e.$
Assume  $\|f_0\|^2_{L^2} < m/2$ and
$m^3/e$ be suitably small.
Then there exists a positive constant $\alpha$
depending on $f_0$, such that
\[
\|f(t)-f_\infty \|_{L^1} \leq C(f_0)  {\e^{-\alpha t/2}},
\]
with  constant $C(f_0)>0.$ Here $f_\infty$ is the Bose distribution with  mass $m$ and energy $e$ defined by \eqref{bose}.
\end{thh}

In comparing the results obtained in \cite{Carrillo} for quantum Fokker-Planck equation, some additional assumptions on the initial data
are needed in Theorem \ref{mainthm}. In fact,
a generalized  logarithmic-Sobolev inequality developed in \cite{carrillo-jungel-el}
for nonlinear diffusion equation was used directly in \cite{Carrillo} to control the entropy production from below
by the relative entropy. In their proof, an auxiliary nonlinear
diffusion equation which has the same entropy and the equilibrium was introduced.
While for the Kac model \eqref{kac},
it is impossible to introduce such auxiliary equation with the same relative entropy or
the equilibrium state, and a compatible entropy production term.
Without using the generalized logarithmic-Sobolev inequality, we follow some ideas used
in \cite{carrillo-jungel-el}, we get the decay rate of entropy production, then the convergence of the solution
of the Kac grazing limit model to its equilibrium. The constraints on the initial data
stated in Theorem \ref{mainthm} will be used to get the decay rate of the entropy
production. We believe that these constrains are only needed
to simplify the proof. While in the numerical simulation part, we don't take into account these
constraints.

The rest of the paper is organised as following. In   Section \ref{s2}, we will give some preliminary estimates. Then the entropy and the entropy equality will be introduced in Section \ref{s3}.
The Bose distribution $f_\infty$ will be given
from the equivalent form of the entropy production.
Then based on a detailed study
on the  entropy production, we get the exponential decay rate.
Finally the theoretical results will be illustrated by  numerical
simulations  in Section \ref{s4}.

\section{Preliminaries}\label{s2}

In this section we will show some a priori weighted $L^2$ estimates on the solution, together with some control on a specific quantity $A_f(t)$.  These estimates will be used
in next section in order to get the decay rate of the entropy production.

Before starting our estimates, we note first that it is not difficult to assert the positivity of the solution $f$ as in \cite{Carrillo}. Let us repeat their arguments quickly here: let $\rho_\eps$ be the Friedrich mollifier, and define the smoothed sign and absolute functions
\[
    \phi_\eps =\rho_{\eps}  \ast  \sign, \quad \Phi_\eps(x) = \int_0^x \phi_\eps(y)\ud y.
\]
Multiplying Kac equation by $\phi_{\eps}(f),$ and integrating it over
$\R$, we get
\begin{multline}
\label{sign}
    \frac {\ud} {\ud t} \int (\Phi_\eps(f) - f) \ud v  = \frac{\ud}{\ud t} \int \Phi_\eps(f)
     \\ =   -A_f(t) \int \phi'_\eps(f) |\partial_v f|^2 \ud v
       - B_f(t) \int v f(1+ f) \phi'_\eps (f) \partial_v f \ud v .
\end{multline}
Note that
\[
    - \int v f \phi'_\eps (f) \partial_v f\ud v
    = -\int v \partial_v (f\phi_\eps (f) -\Phi_\eps(f) ) \ud v ,
\]
and
\[
    - \int v f^2 \phi'_\eps (f)\partial_v f\ud v =
     -\int v\partial_v (f^2 \phi_\eps(f) -  f\Phi_\eps(f)) \ud v
      -\int v \partial_v \left( \int^f \xi^2 \phi'_\eps(\xi) \ud \xi  \right) \ud v.
\]
The first term on the right hand side of \eqref{sign}
is non positive since $\phi'_\eps\geq 0.$
And the second term on the right hand
side of \eqref{sign} vanishes as $\eps \to 0$
from Lebesgure's dominated convergence theorem.
Then  we have by letting $\eps\to 0$ on both side of \eqref{sign},
\[
    \||f(t)|-f(t)\|_{L^1} \leq \||f_0|-f_0\|_{L^1} .
\]
If the initial data $f_0\in L^1$ is non negative a.e on $\R$,
then the solution $f$ (if exists) belongs to $L^1$ and is always non negative a.e on $\R.$

In conclusion, we have shown that if $f$, the smooth solution of Kac's model with initial data $f_0\in L^1(\R)$, is sufficiently decaying, there holds that the $L^1$ norm of $f$ is non-increasing for $t>0.$ Furthermore if $f_0$
is non-negative a.e. in $\R,$  the solution $f(t,v) $ is also non-negative in $\R$ for any $t>0.$

\subsection{Weighted $L^2$ estimates of the solution}
In this paragraph, we are going to show the following
uniform $L^2$ estimates in time.

\begin{prop}
\label{L2-Prop}
The $L^2$ norm of the solution $f$ verifies
\[
    \|f(t)\|^2_{L^2} \leq
     \max \left\{ \|f_0\|^2_{L^2},  \frac{m\Lambda}{3C_N}\sqrt{\frac {m^6}{e^2} + \frac 9{\Lambda^2} \frac{m^3}{e} } \right\}, \quad
     \text{for  } t>0.
\]

Further, assume that $m^3/e$ is sufficient small
and $\|f_0\|_{L^2}^2< m/2,$
then there exists a positive constant $\alpha$
depending on $f_0, \ $ $m$ and $e,$   such that
\begin{equation}
\label{f-L2}
     2 \|f(t)\|_{L^2}^2 - m \leq -\alpha , \quad
     \mbox{for } t>0 .
\end{equation}

Similarly, let $<v>=(1+|v|^2)^{1/2},$ we have the  weighted estimate
\[
    \|<v>f(t)\|_{L^2}^2 \leq
        \max\{ \|<v>f_0\|_{L^2}^2, Y_*(m,e) \},
\]
for $t>0$, with  constant $Y_*(m,e)$ depending on $m$ and $e.$


\end{prop}
\proof
Firstly we multiply the equation by $f$ and then
integrate the resulting equality with respect to $v$
to get
\begin{equation}
\label{L2}
    \frac 1 2 \frac {\ud}{\ud t} \|f(t)\|_{L^2}^2 + A_f(t) \|\partial_v f(t)\|^2_{L^2}
    =m \bigg(\frac 1 2  \|f(t)\|_{L^2}^2 +  \frac 1 3\|f(t)\|_{L^3}^3   \bigg).
\end{equation}
Then we use the interpolation inequality and the Gagliardo-Nirenberg  inequality (cf. (3.27) in \cite{Taylor3}) to
estimate the $L^3$ norm of $f$ as
\begin{eqnarray*}
    \|f(t)\|_{L^3}^3 \leq \|f(t)\|_{L^1} \|f(t)\|_{L^{\infty}}^2 & \leq &  \Lambda \|f(t)\|_{L^1}\|f(t)\|_{L^2}\|\partial_v f(t)\|_{L^2} \\
                            & \leq & \eps  \|\partial_v f(t)\|^2_{L^2} + \frac{\Lambda^2 m^2}{4\eps} \|f(t)\|_{L^2}^2,
\end{eqnarray*}
where $\Lambda$ denotes the constant arising in the Gagliardo-Nirenberg  inequality. Recall that
$A_f(t)=e +\|vf\|_{L^2}^2.$ Choose  $\eps = \frac{3e}{2m}$ and finally we  have
\begin{equation}
\label{L2-1}
    \frac {\ud}{\ud t} \|f(t)\|_{L^2}^2 + e  \|\partial_v f(t)\|^2_{L^2}
    \leq \left( m+ \frac{\Lambda^2 m^4}{9e} \right)  \|f(t)\|_{L^2}^2.
\end{equation}

The Nash inequality (cf. (6) in \cite{Toscani}) in one dimensional case reads
\[
    \|\partial_v f\|_{L^2} \geq C_{  N} \frac {\|f\|^3_{L^2}}{\|f\|^2_{L^1}},
\]
where $C_{  N}$ is a numerical constant.
Using the Nash inequality in \eqref{L2-1} gives
\[
     \frac {\ud}{\ud t} \|f(t)\|_{L^2}^2 \leq -C_2(m,e) \| f(t)\|^6_{L^2}+  C_1(m,e) \|f(t)\|_{L^2}^2,
\]
with $C_1(m,e) = m+ \Lambda^2m^4/(9e),$ $C_2(m,e)= C_{  N}^2e/m^4 . $
The above differential inequality can be solved explicitly in a standard way.
For simplicity let us feel free to omit the dependence of the constants on the mass $m$ and the energy $e$ and denote $X(t) = \|f(t)\|_{L^2}^2$, then we have
\[
    X'(t) \leq -C_2X^3(t) +C_1 X(t),
\]
which can be reduced to
\[
    \left( \frac {\e^{2C_1t}} {X(t)^2} \right)' \geq 2C_2 \e^{2C_1t}.
\]
Integrating the above inequality over $[0,t]$ gives the upper bound of $X^2$ as
\[
    X^2(t) \leq \frac {\e^{2C_1t}} {\frac  1 {X(0)^2}
    + \frac {C_2}{C_1}\bigg( \e^{2C_1t}-1 \bigg)}
    \leq \max \left \{ X(0)^2, \frac {C_1}{C_2} \right\}.
\]
Inserting the expressions of $C_1$ and $C_2$ in the above
inequality, we get the uniform bound for the $L^2$ norm
as stated in Proposition \ref{L2-Prop}.

Next, it is easy to check that
when $m^3/e$ is small enough such that $ \frac {C_1}{C_2} < m/2$, we
get the existence of the positive constant $\alpha$ stated in Proposition \ref{L2-Prop} thus the inequality \eqref{f-L2} holds.

The $L^2$ estimation of  $<v>f$ is similar.
 Let $g=<v>f$ with $<v>= (1+v^2)^{1/2}.$ We multiply the Kac model
\eqref{kac} by $<v>g,$ then integrate it with respect to $v$ over
$\R$ to get
\begin{eqnarray*}
    \frac 1 2 \frac{\ud}{\ud t} \|g(t)\|_{L^2}^2
    & = & A_f(t) \int <v> g \partial_{vv} f \ud v + B_f(t) \int <v>g \partial_v (vf(1+f)) \ud v\\
    & = & A_f(t) \left( \int \frac{v^2}{<v>^4}g^2 \ud v  - \|\partial_v g\|_{L^2}^2 \right) \ud v \\
     & &  + \frac{B_f(t)} 2 \int \frac{g^2}{<v>^2}(1- v^2) \ud v
     + \frac{B_f(t)} 3 \int \frac{g^3}{<v>^3}(1-3v^2) \ud v  .
\end{eqnarray*}
Since $1/<v>\leq 1, \ |v| /<v>  \leq 1,$ we have
\[
    \frac 1 2 \frac {\ud }{\ud t }\|g(t)\|_{L^2}^2
    +A_f(t) \|\partial_v g(t)\|_{L^2}^2 \leq \left( A_f(t)
        + \frac { B_f(t)}2 \right) \|g(t)\|_{L^2}^2
        + \frac { B_f(t)}3 \|g(t)\|_{L^3}^3 .
\]
Similarly we use the interpolation inequality and the
Gagliardo-Nirenberg  inequality to estimate the $L^3$ norm as
\begin{eqnarray*}
    \|g(t)\|_{L^3}^3 & \leq & \Lambda \|g(t)\|_{L^1} \|g(t)\|_{L^2} \|\partial_v
    g(t)\|_{L^2} \\
    & \leq &  \eps \|\partial_v g(t)\|_{L^2}^2 + \frac{\Lambda^2 \|g(t)\|_{L^1}^2}{4\eps}\|
    g(t)\|_{L^2}.
\end{eqnarray*}
We take $\eps = 3e/(2m)$ and use Nash inequality  to get
\[
    \frac {\ud }{\ud t }\|g(t)\|_{L^2}^2
    \leq \left( 2 A_f(t)  + { B_f(t)} +
    \frac{ m^2\Lambda^2 \|g(t)\|^2_{L^1}}{9e} \right) \|g(t)\|_{L^2}^2
    - \frac{eC_N^2}{\|g(t)\|_{L^1}^4}\|g(t)\|^6_{L^2}.
\]
As $<v>\leq 1+ v^2/2, $ we have $ \|g(t)\|_{L^1} \leq m +\frac e 2.$
Note that $A_f(t) \leq  e + \|g(t)\|_{L^2}^2.$ Then we have
\[
    \frac {\ud }{\ud t }\|g(t)\|_{L^2}^2 \leq  \tilde  C_1(m,e) \|g(t)\|_{L^2}^2
    +2\|g(t)\|_{L^2}^4
    -\tilde  C_2(m,e) \|g(t)\|_{L^2}^6,
\]
with
\begin{eqnarray*}
\tilde  C_1(m,e) & = & 2 e + m +
    \frac{ m^2\Lambda^2 (m+\frac e 2)^2 }{9e} , \\
\tilde  C_2(m,e) & = & \frac{eC_N^2}{(m+\frac e2)^4}.
\end{eqnarray*}
Let $Y(t) = \|g(t)\|_{L^2}^2 $ and the differential inequality can
 be written as
\[
  Y'(t) \leq G(Y(t)) = \tilde  C_1(m,e)  Y(t)+ 2Y(t)^2 -\tilde
  C_2(m,e)Y^3(t).
\]
As $G(Y)$ has a unique positive zero point
$$Y_*(m,e)= \frac{1+\sqrt{1+ \tilde  C_1(m,e)\tilde  C_2(m,e)}}{\tilde  C_2(m,e)}>0,$$
and $G(Y)$ is
positive over $]0,Y_*[$ and negative on $]Y_* , + \infty[.$ Then we
get the global existence of $Y(t)$ which will take values between
the initial value $Y(0)$  and the equilibrium point $Y_*.$
In conclusion, we have
\[
    \|g(t)\|_{L^2}^2 \leq \max\{\|g_0\|_{L^2}^2, Y_*\},
\]
with $g_0= <v>f_0.$

\begin{remarque}

\begin{enumerate}
\item Note that we have the following interpolation inequality
\[
    \|f(t)\|^2_{L^2} \geq \frac m {2^{7/2}} \left( \frac {m^3} e \right)^{1/2},
\]
which follows classically by optimizing w.r.t. $R>0$ the following inequality
\[
    m = \int_{|v|\leq R} f\ud v + \int_{|v|\geq R} f \ud v
    \leq  \|f\|_{L^2} \sqrt{2R} + \frac e {R^2}.
\]

Therefore an $L^2$ control on $f$ implies a control of $\displaystyle {m^{5/2}\over e^{1/2}}$.
\item We have only shown weighted $L^2$ estimations of $f$ but similar estimates on higher derivatives also hold true. It is important to note that smoothness is not required for estimating convergence to equilibrium.
\end{enumerate}
\end{remarque}

\subsection{Estimate for the quantity $A_f(t)$ }
In this paragraph we study the time derivative of the quantity $A_f(t)= e + \|vf\|_{L^2}^2.$
\begin{prop}\label{A(t)}
There holds
\[
    \frac{A_f'(t)}{A_f(t)} \leq 2\|f(t)\|_{L^2}^2.
\]
\end{prop}
\proof
Let
$g(t,v) = vf(t,v),$ which verifies
\[
    \partial_t g = A_f(t) \partial_{vv} g-2 A_f(t) \partial_{v}f
    + B_f(t) \partial_{v}(vg+g^2) - B_f(t)g(1+f).
\]
Then multiply the equality by $g$ and integrate it with respect to $v$. Finally we get the
$L^2$ equality as
\begin{multline*}
    \frac 1 2 \frac{\ud}{\ud t} \|g(t)\|^2_{L^2}
    + A_f(t) \|\partial_v g\|^2_{L^2} +\frac {B_f(t)}2 \|g(t)\|^2_{L^2} \\
    = A_f(t) \|f(t)\|^2_{L^2} - B_t(f) \int f(t,v)g^2(t,v) \ud v.
\end{multline*}
As $A_f(t)=e + \|g(t)\|^2_{L^2},$ we have
\begin{equation*}
A_f'(t)  +A_f(t) (B_f(t)- 2 \|f(t)\|^2_{L^2}) \leq B_f(t) e  \leq  B_f(t) \cdot  A_f(t).
\end{equation*}
Using the definition of $A_f(t)$ we complete the proof of Proposition \ref{A(t)}.

\section{Relative entropy method and decay to equilibrium}
\label{s3}

In this section, we will prove Theorem \ref{mainthm} by the
relative entropy method. Firstly we will introduce
the entropy, the entropy production and the equilibrium to the Kac
model \eqref{kac}.  Secondly we will show the decay rate of the
entropy production. Finally the decay rate of the solution for
the Kac model to the equilibrium can be derived.
\subsection{Entropy and equilibrium}
Let $\gamma(f) =  f\log f - (1+f)\log (1+f).$ Note that
$\gamma'(f) = \log \frac f {1+f}.$
The entropy $H(f)$
\[
    H(f) =\int \gamma (f) \ud v
\]
verifies the  entropy equality
\begin{equation}\label{ent}
    \frac{\ud}{\ud t} H(f) =
        - \int \left(A_f(t) \frac {|\partial_v f|^2}{f(1+f)} + B_f(t)v \partial_v f \right): = -D(f).
\end{equation}
The entropy production $D(f)$ can be written in some other forms.
For example since $B_f(t) = \int_v f = -\int_v v\partial_v f,$
$D(f)$ can be written as
\begin{eqnarray}
D(f)& =  & \int \left( A_f(t) \frac {|\partial_v f|^2}{f(1+f)} + 2B_f(t) v\partial_v f +B(t)^2  \right )\ud v \nonumber \\
    & =  & \frac 1 {A_f(t)} \int f(1+f) \left| A_f(t) \partial_v \gamma'(f)  + B_f(t) v  \right |^2\ud v. \label{EntroProd-1}
\end{eqnarray}
From \eqref{EntroProd-1} we get that $D(f)(t) \geq 0.$

Furthermore, we can
use the expression of $A_f(t)$ into the entropy
production $D(f)$ and write it in a symmetric form as
\begin{eqnarray*}
    D(f) & =& \int_v \!\!\int_{v_*} \left(  v_*^2 f_*(1+f_*)  \frac {|\partial_v f|^2}{f(1+f)} -  v\partial_v f v_* \partial_{v_*}f_* \right) \ud v \ud v_* \\
    & = & \frac 1 2 \int_v \!\!\int_{v_*}   f(1+f)f_*(1+f_*) \left| v  \partial_{v_*} \gamma'(f_*) -  v_* \partial_v \gamma'(f)  \right|^2 \ud v \ud v_* ,
\end{eqnarray*}
with $f_*=f(v_*).$ From the equality
\[
     v_* \partial_v \gamma'(f) = v \partial_{v_*} \gamma'(f_*) ,
\]
the equilibrium $f_{\infty}$ is
\begin{equation}
\label{bose}
    f_\infty(v) = \frac 1 {\exp(\lambda_1 v^2 - \lambda_2 )-1},
\end{equation}
where the constants $\lambda_1 >0$ and $\lambda_2$ will be determined by the initial data. Note that the equilibrium $f_\infty$ defined
above is the so-called Bose distribution function.

\begin{remarque}We can show that
\begin{equation}\label{lambda1}
    \lambda_1 = \frac {B_{f_\infty}}{2A_{f_\infty}}.
\end{equation}
In fact, we have
\begin{eqnarray*}
    B_{f_\infty} &=&    \frac 1 {\sqrt{\lambda_1}} \int  \frac {e^{\lambda_2}} {e^{  v^2}- e^{\lambda_2} } \ud v \\
    A_{f_\infty} &=&
    \int \frac{v^2 \e^{\lambda_1 v^2} e^{\lambda_2}}{(e^{\lambda_1 v^2}- e^{\lambda_2 })^2} \ud v
    = - \partial_{\lambda_1} \int \frac {e^{\lambda_2}} {e^{\lambda_1 v^2}- e^{\lambda_2} } \ud v.
\end{eqnarray*}
Using the expressions of $B_{f_\infty},$ we have
\[
    A_{f_\infty} =  - \partial_{\lambda_1} B_{f_\infty}  =   \frac 1 {2 \lambda_1^{3/2}} \int  \frac {e^{\lambda_2}} {e^{  v^2}- e^{\lambda_2} } \ud v = \frac 1 {2\lambda_1} B_{f_\infty} .
\]
Note that \eqref{lambda1} can also be obtained from the entropy production in the form \eqref{EntroProd-1}.
\end{remarque}

In conclusion, we have the following lemma
\begin{lemm}
The equilibrium $f_\infty$ minimizes
\[
    \left\{H(f): f(v) \mbox{ is positive, }\int f(v) \ud v =m, \ \int v^2 f(v) \ud v =e \right\}
\]
with $m$ and $e$ fixed. As $\gamma$ is convex, this minimizer function is unique.
Moreover, given any solution $f(t,v)$ to the Kac model \eqref{kac} with initial data
$f_0$ of mass $m$ and energy $e,$ we have
\[
    H(f_\infty) \leq H(f)(t) \leq H(f_0), \quad t>0,
\]
and
\[
    \lim_{t\to\infty} H(f)(t) = H(f_\infty).
\]
\end{lemm}

Before ending this paragraph, we introduce the relative entropy $H(f|f_\infty)$ as
\[
    H(f|f_\infty) = H(f)-H(f_\infty) = \int \left[\gamma(f)- \gamma(f_\infty) -\gamma'(f_\infty)(f-f_\infty)   \right] \ud v,
\]
where we used the conservations of mass and energy for the last equality.

\subsection{Decay rate of the entropy production and the relative entropy}

To get the decay rate of the entropy production, we shall study the
time derivative of $D(f)$. To simplify notations we denote
$\xi = A_f(t) \partial_v \gamma'(f) +B_f(t) v. $
Hence the Kac equation and the entropy production $D(f)$
can be written as
\[
    \partial_t f =\partial_v [ f(1+f) \xi], \quad
    D(f) = \frac 1 {A_f(t)} \int f(1+f) \xi^2 .
\]
Then we have
\begin{multline}
\label{EntroProd-deri}
\frac \ud {\ud t}  D(f) = -\frac {A_f'(t)}{A_f^2(t)} \int f(1+f)\xi^2  + \frac 1 {A_f(t)}\int (1+2f)\xi^2 \partial_t f
+ \frac 2 {A_f(t)} \int f(1+f) \xi \partial_t \xi \\
        := \mbox{I} +\mbox{II} + \mbox{III} .
\end{multline}
Next we will calculate these three integrals. Firstly the integral $\mbox{I}$ can be written as
\[
    \mbox{I} = -\frac {A_f'(t)}{A_f(t)} D(f).
\]
Then the second integral $\mbox{II}$ can be calculated as
\[
    \mbox{II} = \frac {1}{A_f(t)} \int (1+2f) \xi^2 \partial_v [f(1+f)\xi]
    = - \frac 2 {A_f(t)} \int \left(f^3 +\frac 3 2 f^2 +f\right )\xi^2 \partial_v \xi \ud v.
\]
We denote by $\phi (f) = f^3 +\frac 3 2 f^2 +f.$ Then using the expression of
$\xi= A_f(t)\partial_v \gamma'(f) +B_f(t) v,$ we can rewrite $\mbox{II}$ as
\begin{multline}
    \mbox{II} = -\frac{B_f(t)}{A_f(t)}\int \phi(f) \xi^2 \ud v + 2 \int \left( \frac {\phi'(f)}{f(1+f)} -2 \frac {\phi(f) (1+2f)}{f^2 (1+f)^2} \right)|\xi\partial_v f|^2 \\
        + 4 \int \frac {\phi(f) }{f^2 (1+f)^2}  \xi\partial_v f \partial_v [f(1+f)\xi].
\end{multline}
Finally, from the conservation of mass, $B_f(t)=m,$
we get
\[
    \mbox{III} = 2\frac{A_f'(t)}{A_f(t)} \int \xi \partial_v f  +  2 \int f(1+f) \xi \partial^2_{tv} \gamma'(f)
    =  2\frac{A'_f(t)}{A_f(t)} D(f) - 2 \int \frac 1{f(1+f)}|\partial_v [f(1+f)\xi]|^2.
\]

In summarize, we get the derivative of the entropy dissipation
\[
    \frac \ud {\ud t} D(f) =  -\frac{B_f(t)}{A_f(t)}\int \phi(f) \xi^2 \ud v +  \frac{A'_f(t)}{A_f(t)} D(f) - \{ \cdots\},
\]
where $\{\cdots\}$ denotes some positive terms. As $\phi(f) >f(1+f),$ then we get
\[
    \frac \ud {\ud t} D(f) \leq  \left( \frac { A'_f(t)} {A_f(t)} - B_f(t) \right) D(f) .
\]

We use first Proposition \ref{A(t)} and then Proposition \ref{L2-Prop} to get
\[
    \frac{A'_f(t)}{A_f(t)} -B_f(t)  
     \leq  2\|f(t)\|_{L^2}^2  -B_f(t)\leq -\alpha.
\]
In the last inequality we used the smallness assumptions as
in Proposition \ref{L2-Prop}.
Immediately we derive the following decay rate about the
entropy production
\[
    D(f)(t) \leq D(f_0) \e^{-\alpha t}.
\]

We use the decay rate of the entropy production $D(f)$ 
in the entropy equality \eqref{ent} to get 
\[
    \frac{\ud }{\ud t}H(f|f_\infty) \geq -D(f_0)\e^{-\alpha t}.
\]
Then integrating the above inequality over $]t_1,t_2[$ gives
\[
H(f|f_\infty)(t_2) -H(f|f_\infty)(t_1)
\geq  D(f_0) \frac {\e^{-\alpha t_2}-\e^{-\alpha t_1}}{\alpha}.
\]
Let $t_2\to +\infty, $ and  as
  $\lim_{t\to \infty} H(f|f_\infty)(t)=0,$
finally  we get the decay rate of the relative entropy
\[
    H(f|f_\infty)(t) \leq D(f_0) \frac {\e^{-\alpha t}}{\alpha}.
\]

\subsection{Decay rate of the $L^1$ norm and the proof of Theorem \ref{mainthm}}

Next, we show the $L^1$ decay rate of the solution to
the equilibrium.
Observe that there exists  a function $y(t,v)$ which takes values between $f(t,v) $ and $f_\infty(v)$
such that the relative entropy can be written as
\begin{eqnarray*}
    H(f|f_\infty) = H(f)-H(f_\infty)
    & =& \int (\gamma(f) -\gamma(f_\infty) -\gamma'(f_\infty) (f-f_\infty ))\ud v \\
    & =& \int \gamma''(y(t,v))(f-f_\infty)^2 \ud v ,
\end{eqnarray*}
Remark that we used
the property of mass and energy conservations and the Taylor formula in the last
two equalities.

As \cite{Carrillo},  using the Cauchy Schwartz inequality, we have
\begin{eqnarray*}
    \|f(t)- f_\infty\|^2_{L^1(f<f_\infty)}
    & \leq & \int_{\{f<f_\infty\}} \frac 1 {\gamma''(y(t,v))}
       \int_{\{f<f_\infty\}} {\gamma''(y(t,v))} (f-f_\infty)^2  \\
    & \leq & \int f_\infty (1+f_\infty)   \int  {\gamma''(y(t,v))} (f-f_\infty)^2   \\
    & \leq & C\  H(f|f_\infty)   \leq C\ D(f_0) \frac {\e^{-\alpha t}}{\alpha}.
\end{eqnarray*}
Hence using the mass conservation
we get the following desired result
\[
     \|f(t)- f_\infty\| _{L^1(\R)} = 2  \|f(t)- f_\infty\|_{L^1(f<f_\infty)}
     \leq C (f_0)  {\e^{-\alpha\ t/2}}.
\]
Hence the proof to Theorem \ref{mainthm} is completed.

\section{Convergence towards equilibria: numerical simulations}
\label{s4}

Since we have shown the solution goes exponentially fast to the Bose equilibrium distribution, we will do some numerical simulations in this section, to show the equilibrium distributions for different initial states, and the exponential decay of the entropies.

We recall first the Bose distribution
\[
    f_\infty =\frac 1 {\exp(\lambda_1 v^2 -\lambda_2)-1}.
\]

The numerical simulations are carried out by different initial conditions. The first example shows if the initial data is concentrated near the center, it will evolve to Bose distribution, with entropy decaying exponentially to some final state.

\begin{figure}[!h]
\begin{minipage}[t]{0.5\linewidth}
\centering

\subfigure{
\includegraphics[width=3.2in]{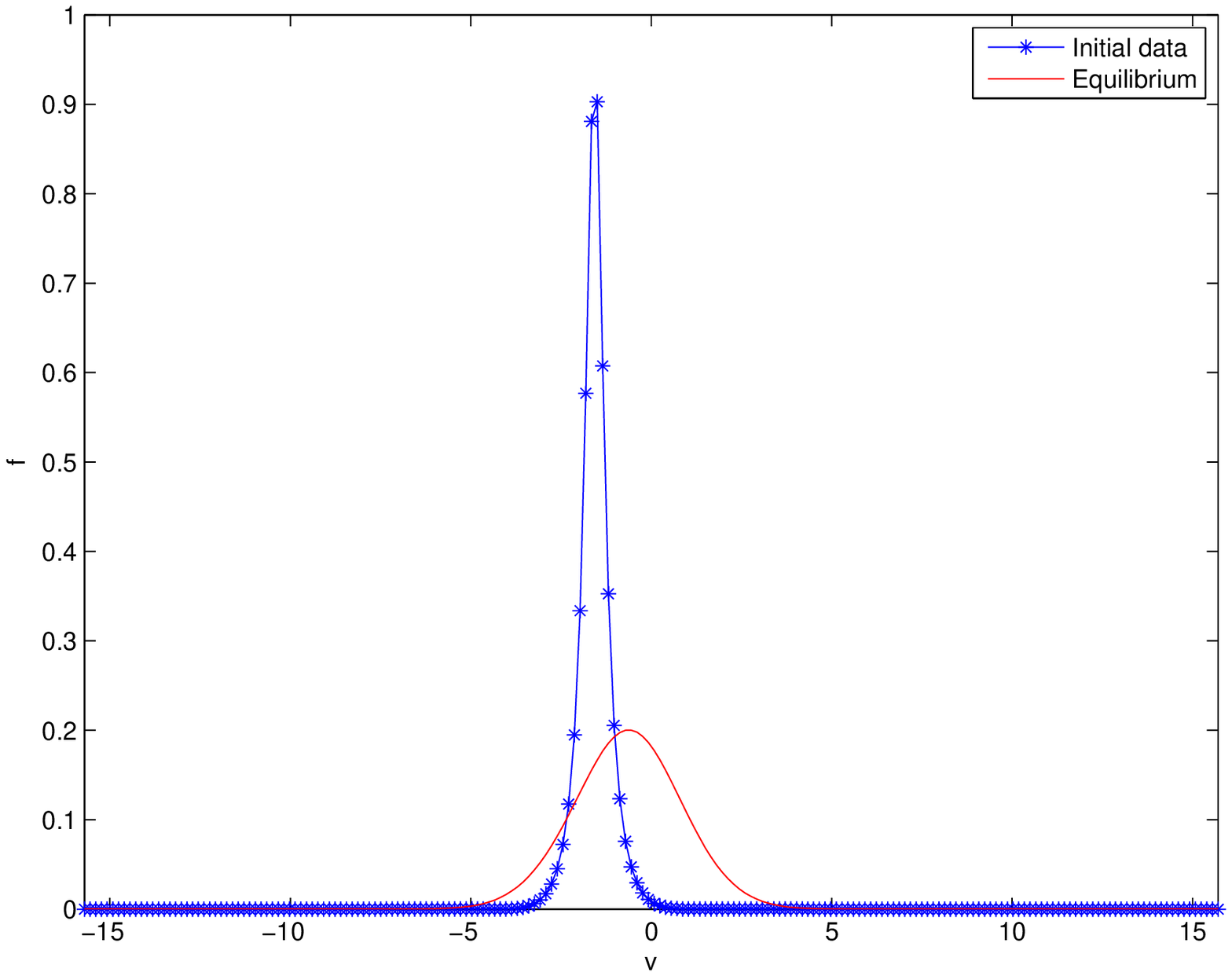}}

\end{minipage}%
\begin{minipage}[t]{0.5\linewidth}
\centering

\subfigure{
\includegraphics[width=3.2in]{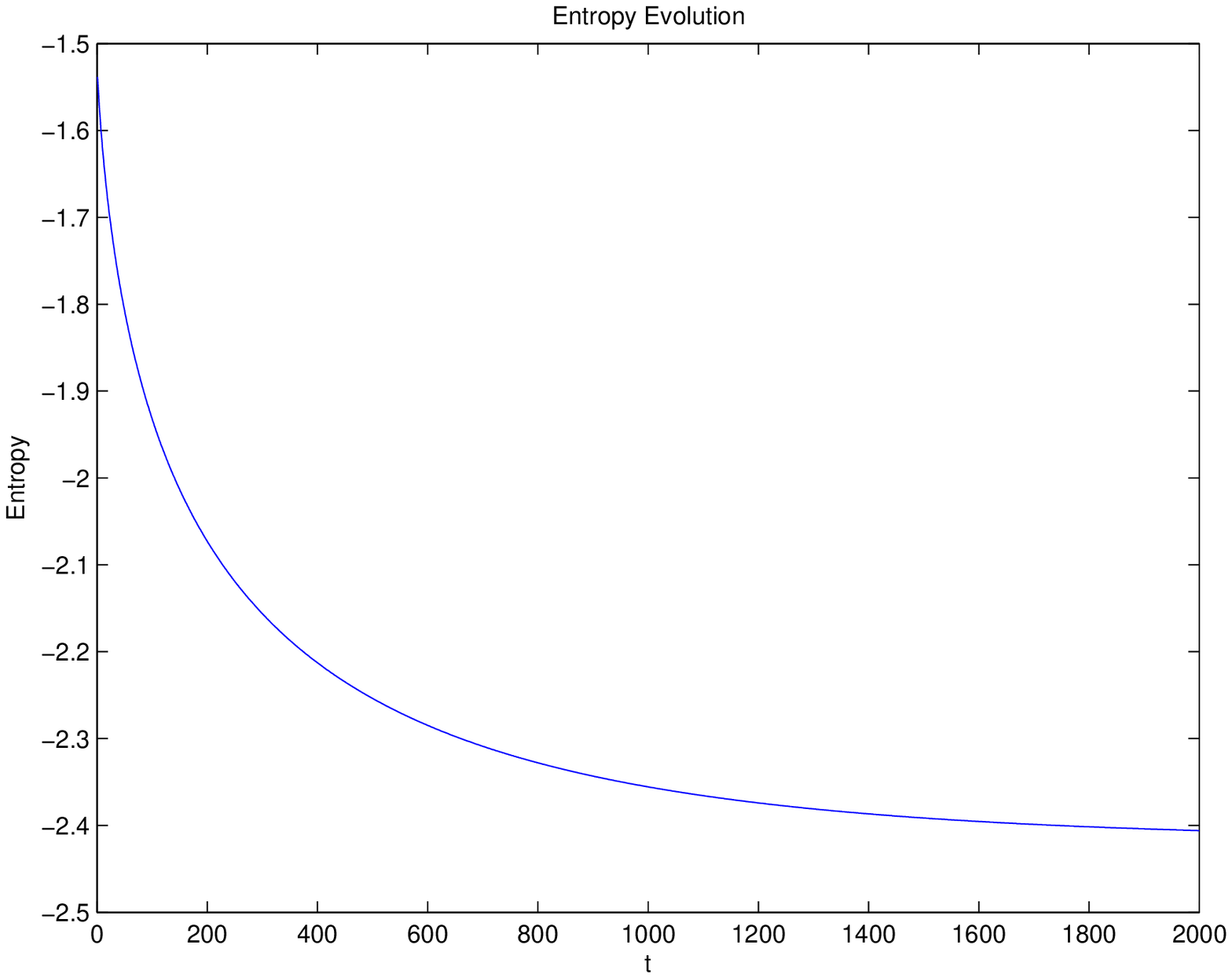}}

\end{minipage}

\caption{Initial data $f_0$ in (\ref{init1}).}
\label{fig1}
\end{figure}

\noindent {\bf Example 1.}  Consider initial data
 \begin{equation}\label{init1}
 f_0=\dfrac{0.1}{e^{(v-\pi /2)^2+0.1}-1}.
 \end{equation}
The equilibrium distribution and evolution of entropy are shown in Figure \ref{fig1}.

This property is true without smallness assumption on the initial data. In Example 2, we take 10 times the initial data as in Example 1 and observe also the exponential decay of the entropy, with different time scale used in the simulation.

\begin{figure}[!h]
\begin{minipage}[t]{0.5\linewidth}
\centering

\subfigure{
\includegraphics[width=3.2in]{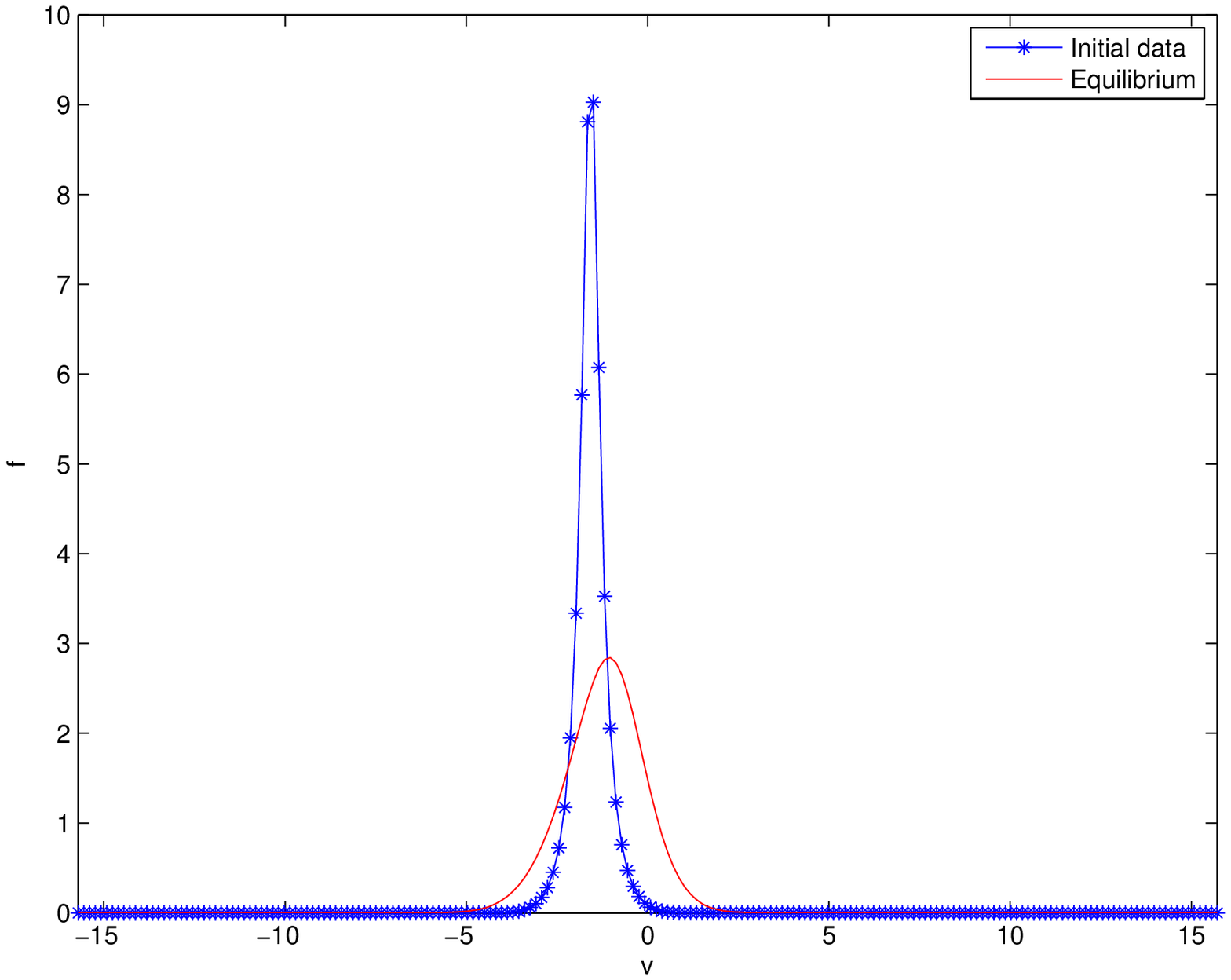}}

\end{minipage}%
\begin{minipage}[t]{0.5\linewidth}
\centering

\subfigure{
\includegraphics[width=3.2in]{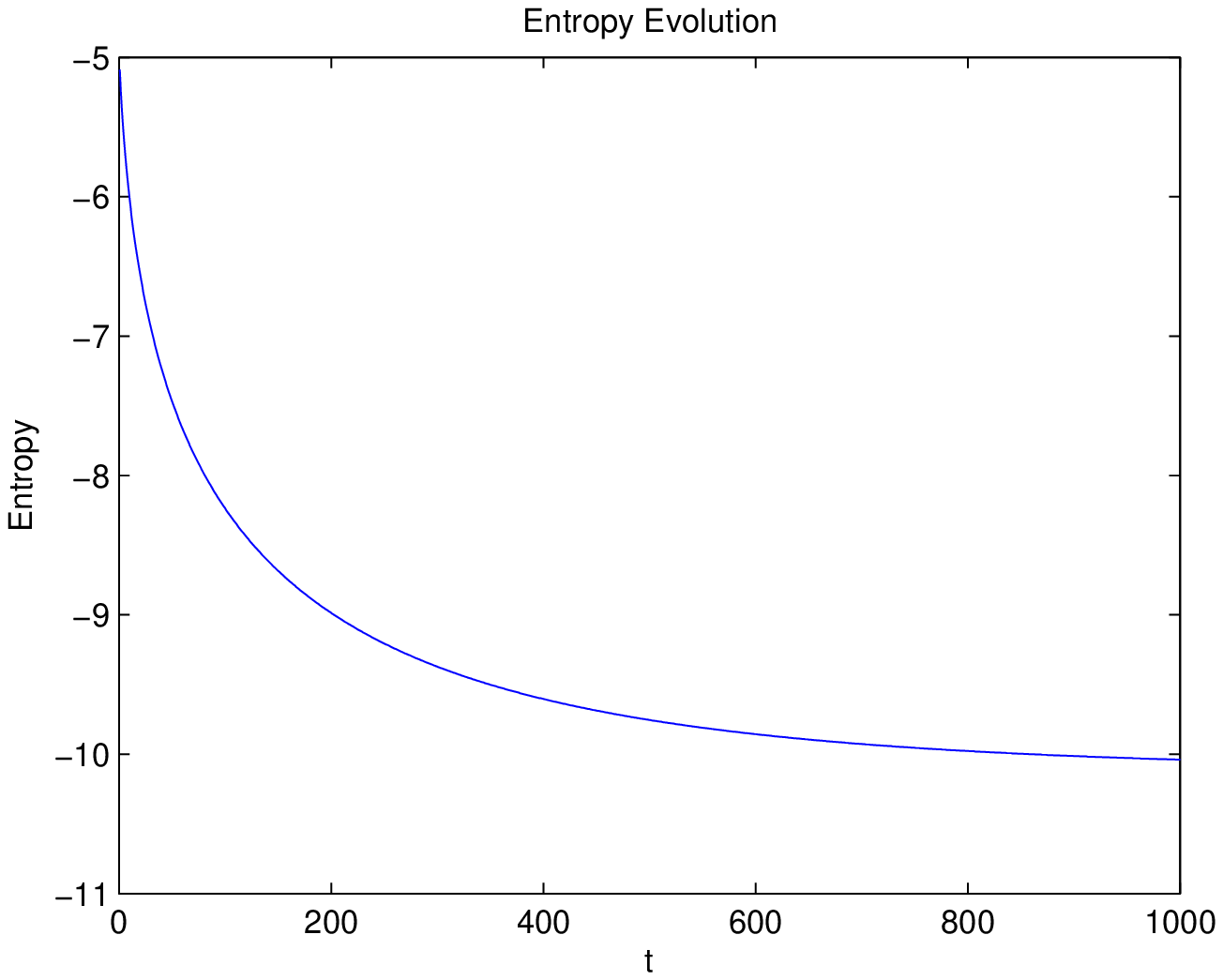}}

\end{minipage}

\caption{Initial data $f_0$ in (\ref{init5}).}
\label{fig2}
\end{figure}

\noindent {\bf Example 2.}  Consider initial data
  \begin{equation}\label{init5}
 f_0=\dfrac{1}{e^{(v-\pi /2)^2+0.1}-1}.
 \end{equation}
The equilibrium distribution and evolution of entropy are shown in Figure \ref{fig2}.
\vspace{2mm}

 As we know the Bose distribution behaves like Gaussian when $|v|$ is big. The next example shows the evolution of a Gaussian to Bose distribution.

\begin{figure}[!h]
\begin{minipage}[t]{0.5\linewidth}
\centering

\subfigure{
\includegraphics[width=3.2in]{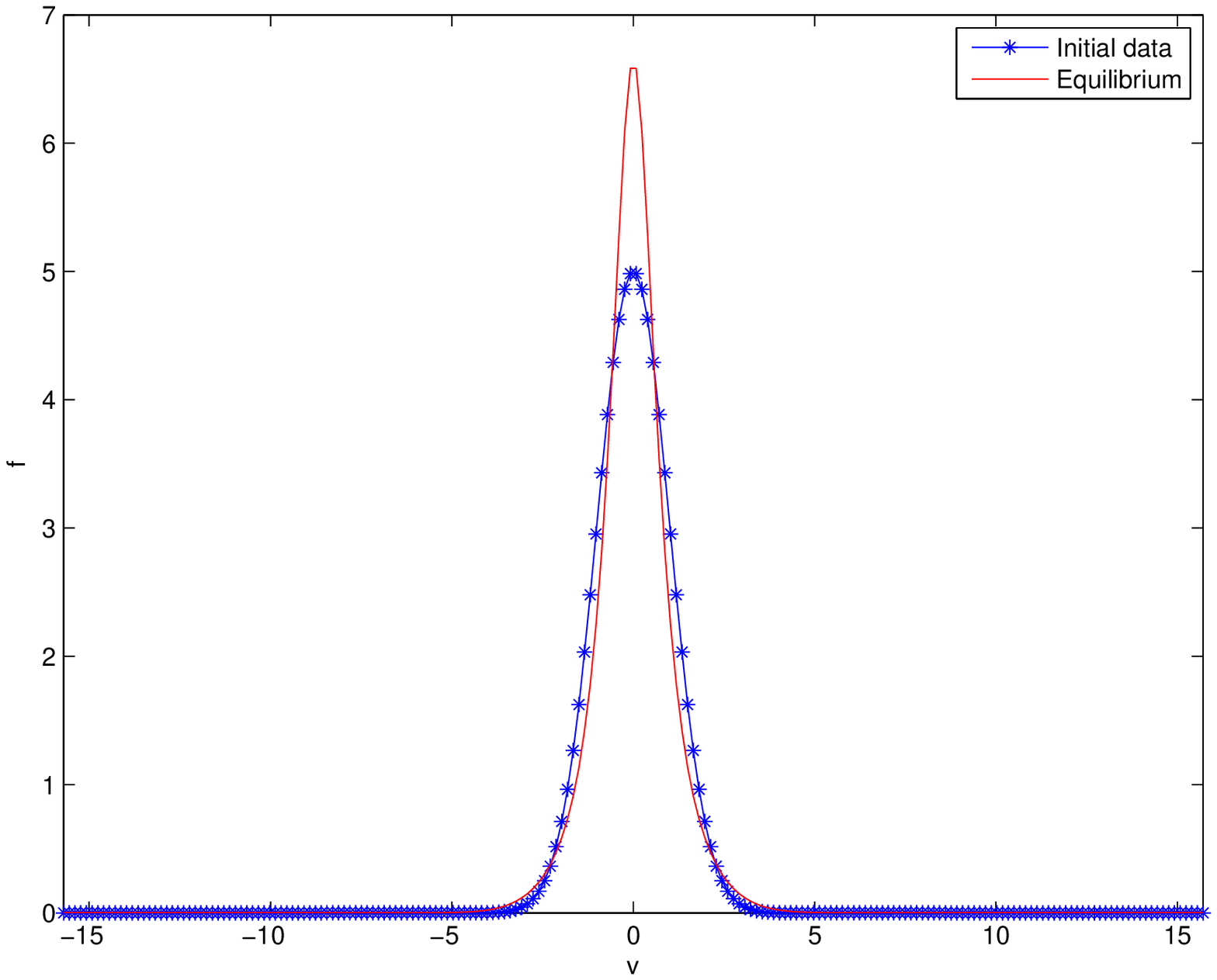}}

\end{minipage}%
\begin{minipage}[t]{0.5\linewidth}
\centering

\subfigure{
\includegraphics[width=3.2in]{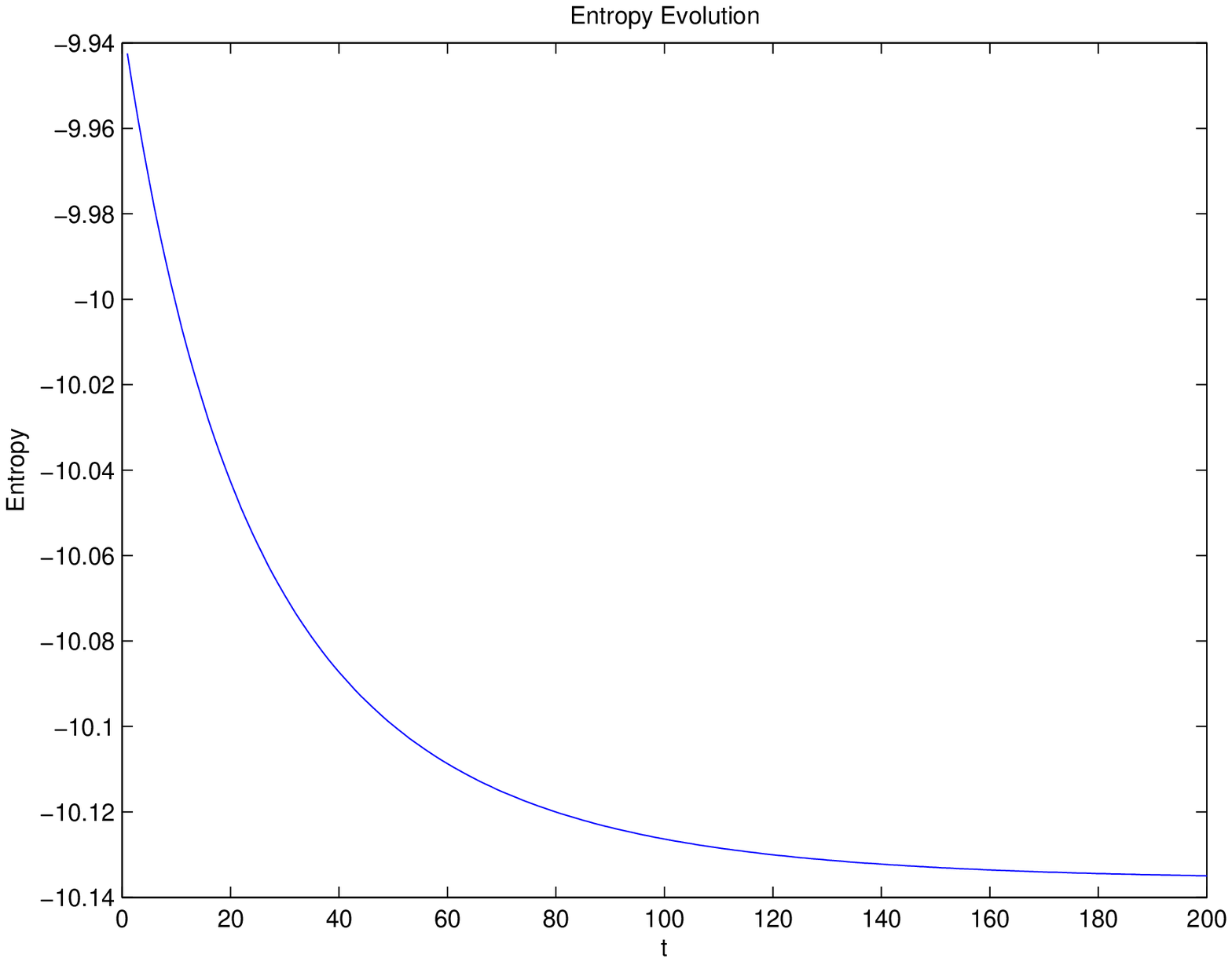}}

\end{minipage}

\caption{Initial data $f_0$ in (\ref{init2}).}
\label{fig3}
\end{figure}

\noindent {\bf Example 3.}  Consider initial data
  \begin{equation}\label{init2}
 f_0=5*e^{-v^2/2}.
 \end{equation}
The equilibrium distribution and evolution of entropy are shown in Figure \ref{fig3}.

The comparison shows the Bose distribution is more singular near $|v|=0$ , but bahaves like a Gaussian for $|v|$ big.

\begin{figure}[!h]
\begin{minipage}[t]{0.5\linewidth}
\centering

\subfigure{
\includegraphics[width=3.2in]{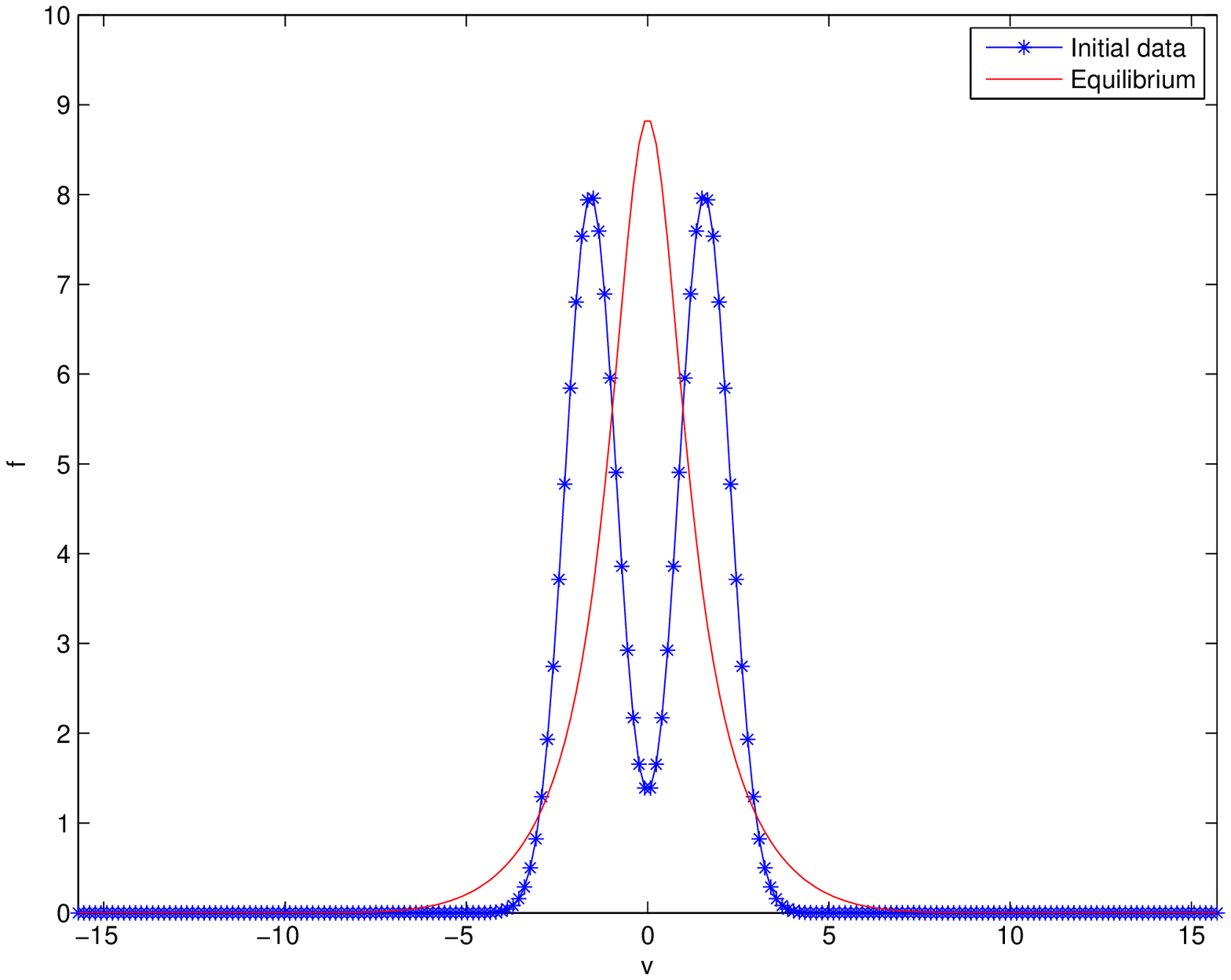}}

\end{minipage}%
\begin{minipage}[t]{0.5\linewidth}
\centering

\subfigure{
\includegraphics[width=3.2in]{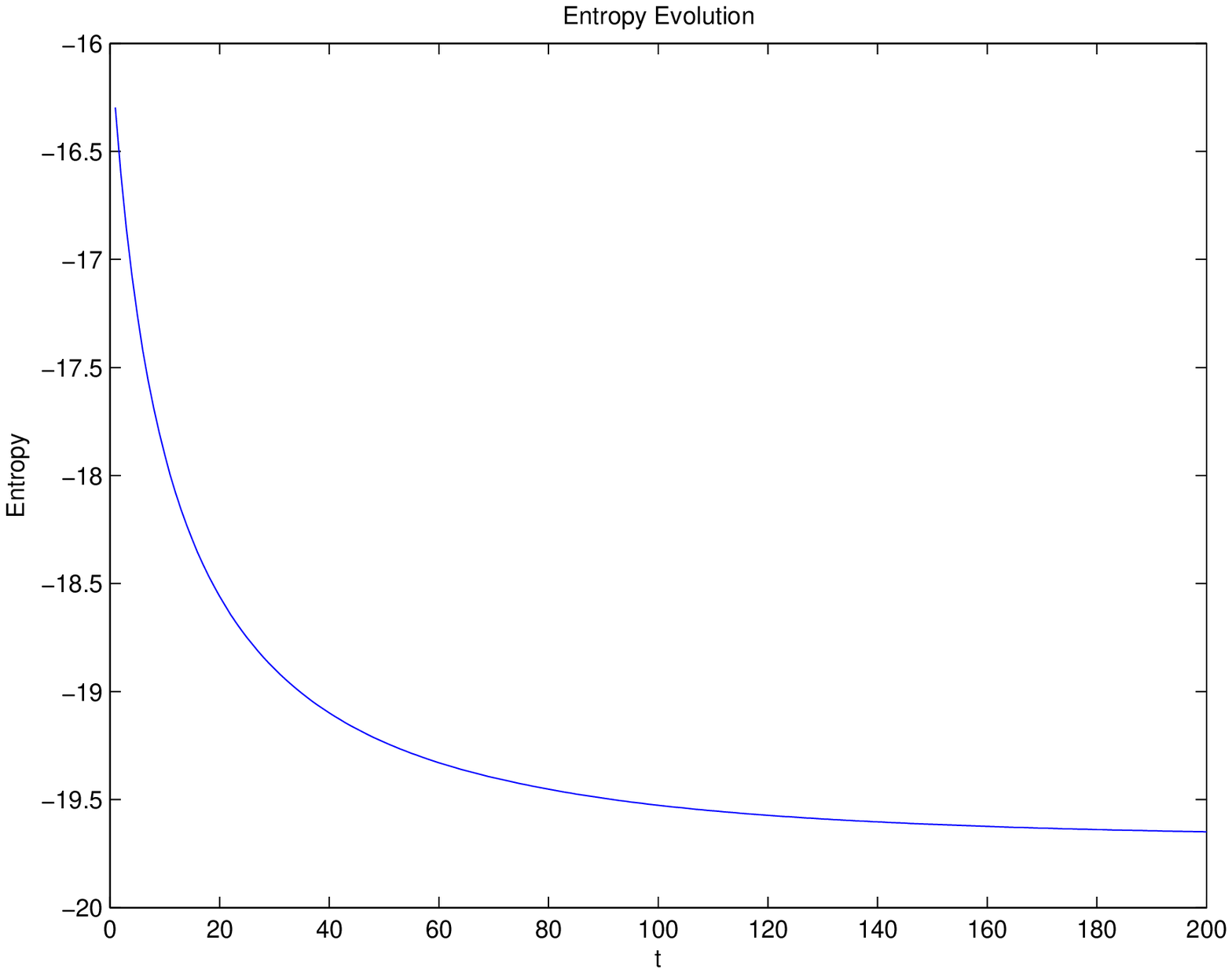}}

\end{minipage}

\caption{Initial data $f_0$ in (\ref{init3}).}
\label{fig4}
\end{figure}

 \noindent {\bf Example 4.}  Consider initial data
 \begin{equation}\label{init3}
 f_0=8*[e^{-(v+\pi /2)^2}+e^{-(v-\pi /2)^2}].
 \end{equation}
The equilibrium distribution and evolution of entropy are shown in Figure \ref{fig4}.

This example shows the evolution of summation of two Gaussians. We will show a more general case in next example.

\begin{figure}[!h]
\begin{minipage}[t]{0.5\linewidth}
\centering

\subfigure{
\includegraphics[width=3.2in]{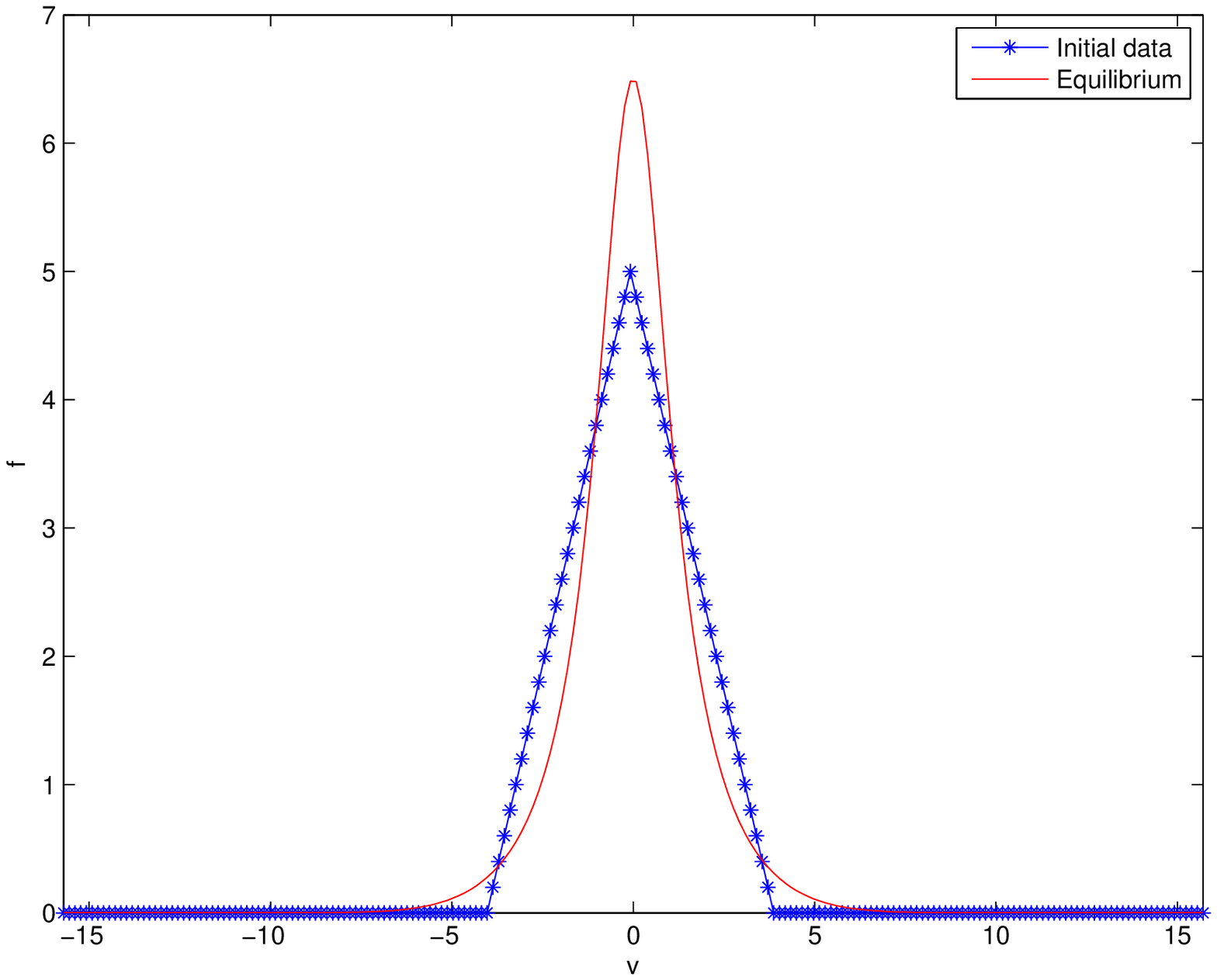}}

\end{minipage}%
\begin{minipage}[t]{0.5\linewidth}
\centering

\subfigure{
\includegraphics[width=3.2in]{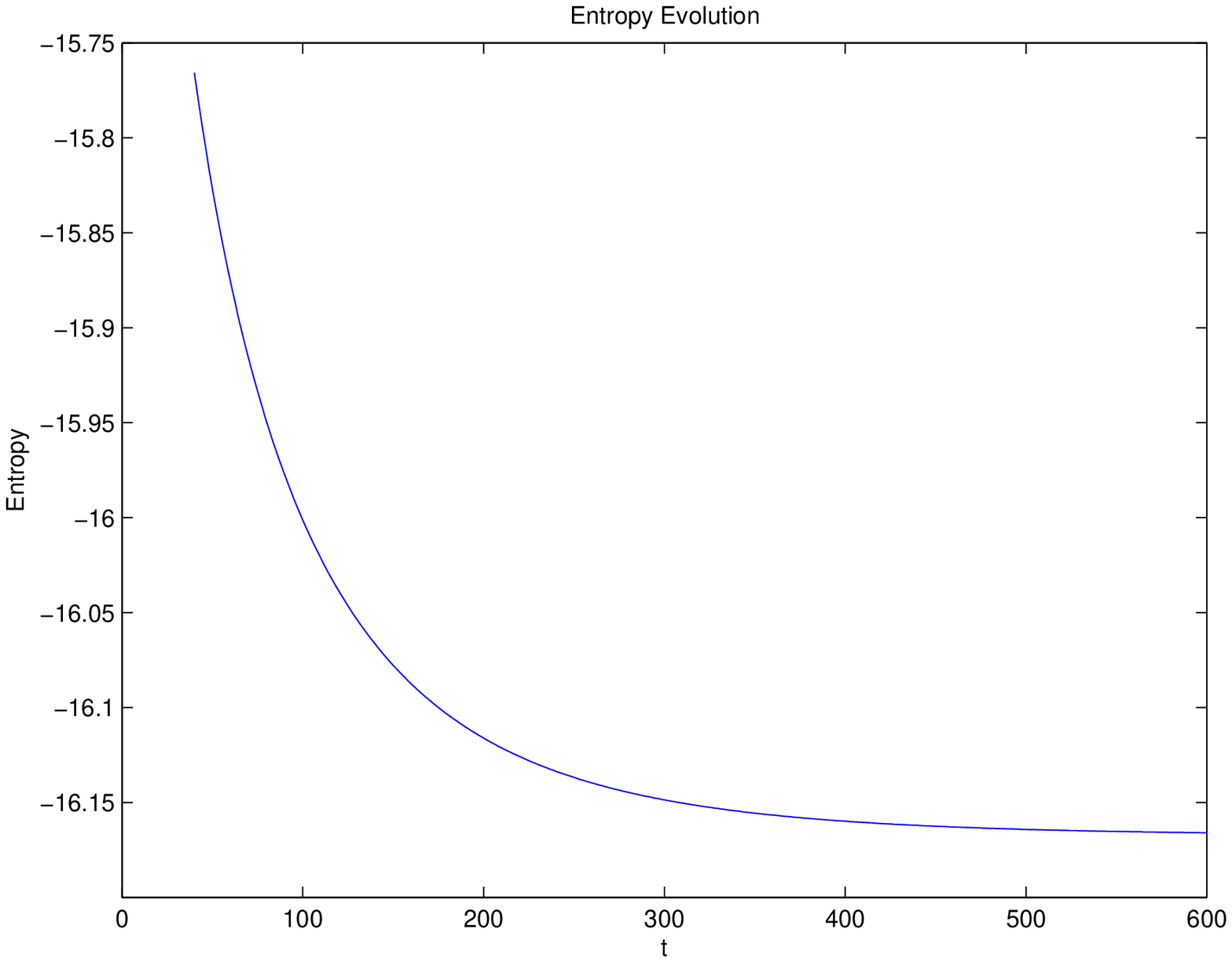}}

\end{minipage}

\caption{Initial data $f_0$ in (\ref{init4}).}
\label{fig5}
\end{figure}

\noindent {\bf Example 5.}  Consider initial data
 \begin{equation}\label{init4}
 f_0=\left\{
 \begin{array}{ll}
  \dfrac{5}{2}+\dfrac{2}{\pi}v \text{~~for~~} v\in [-\dfrac{5\pi}{4}, 0],  \vspace{1.5mm}\\
  \dfrac{5}{2}-\dfrac{2}{\pi}v \text{~~for~~} v\in [0, \dfrac{5\pi}{4}],  \vspace{1.5mm}\\
  0 \text{~~for~~others}.
\end{array}
  \right.
 \end{equation}
The equilibrium distribution and evolution of entropy are shown in Figure \ref{fig5}.

Note that this initial data is not in space $L\log L$, since $f_0=0$ in some interval thus the entropy is $\infty$ at first several time steps.
 This more general case also shows the exponential decay of entropy.

All the above numerical results showed the quick convergence towards to equilibria, especially, we can see the exponential evolution of entropy. The numerical results further elaborate our main result stated in Theorem \ref{mainthm}.

\section*{Acknowledgement}
This work was supported by the Fundamental Research Funds for the Central Universities and National Natural Science Foundation of China (Nos.11171211, 11171212, 11201116), China Postdoctoral Science Foundation,
 a starting grant from Shanghai Jiao Tong University, 
 together with Shanghai Rising Star Program (12QA1401600).

\bibliographystyle{plain}
\bibliography{Kacbib}
\end{document}